\documentstyle[12pt]{article}
\newtheorem{thm}{Theorem}

\newtheorem{lem}{Lemma}
\newtheorem{prop}{Proposition}

\newcommand{\Ho}{H\"{o}lder}
\newcommand{\T}{Teichm\"{u}ller}

\newcommand{\WonesG}{${\cal W}_{1}^{s,G}$}
\newcommand{\WoneuG}{${\cal W}_{1}^{u,G}$}
\newcommand{\WtwosG}{${\cal W}_{2}^{s,G}$}
\newcommand{\WtwouG}{${\cal W}_{2}^{u,G}$}

\begin{document}

\begin{center}
{\bf The \T \ space of the standard action of $SL(2,Z)$ 
on $T^2$ is trivial}\\
\bigskip
\bigskip
Elise E. Cawley\\
Department of Mathematics\\
University of California at Berkeley\\
Berkeley, CA 94720\\
\bigskip
February 12, 1991
\end{center}
\bigskip
\bigskip
\bigskip

\noindent {\bf 1} \ \ \ {\bf Introduction}

\bigskip

\noindent The group $SL(n,{\bf Z})$ acts linearly on ${\bf R}^n$, 
preserving
the integer lattice ${\bf Z}^{n} \subset {\bf R}^{n}$.  The induced
(left) action  on the n-torus 
${\bf T}^{n} = {\bf R}^{n}/{\bf Z}^{n}$ will be referred to as the
``standard action''.
It has recently been shown that the standard action of 
$SL(n,{\bf Z})$ on ${\bf T}^n$, for $n \geq 3$, is 
both topologically and smoothly rigid. \cite{Hu}, \cite{KL}, \cite{HKLZ}.
That is,
nearby actions in the space of representations of $SL(n,{\bf Z})$ 
into ${\rm Diff}^{+}({\bf T}^{n})$ are smoothly conjugate to the standard
action.  In fact, this rigidity persists for the standard action of
a subgroup of finite index.
On the other hand, while the ${\bf Z}$ action on ${\bf T}^{n}$ defined
by a single hyperbolic element of $SL(n,{\bf Z})$ is topologically rigid,
an infinite dimensional space of smooth conjugacy classes occur in 
a neighborhood of the linear action. \cite{C}, \cite{L}, \cite{MM1}, \cite{MM2} 
The standard action of
$SL(2, {\bf Z})$ on ${\bf T}^2$ forms an intermediate case, with different rigidity 
properties from
either extreme.  One can construct continuous deformations of 
the standard action to obtain  an (arbritrarily near) 
action to which it is not topologically conjugate.
\cite{Hu} \cite{Ta}
The purpose of the present paper is to show that
if a nearby action, or
more generally, an action with some mild Anosov properties, is conjugate
to the standard action of $SL(2, {\bf Z})$ on ${\bf T}^2$ by a homeomorphism $h$, 
then $h$ is smooth. In fact, it will be shown that this rigidity holds for 
any non-cyclic
subgroup of $SL(2, {\bf Z})$. 

A smooth action of a group $\Gamma$ on a manifold $M$ is called 
Anosov if there is a least one element $\gamma \in \Gamma$ which
acts by an Anosov diffeomorphism.  We define the $C^{r}$ \T \ space of 
an Anosov action $F:\Gamma \times M \rightarrow M$
to be the space of ``marked'' $C^r$ smooth structures 
preserved by the underlying topological dynamics of the action, and such that
Anosov group elements for the action $F$ define Anosov diffeomorphisms 
with respect to 
the new smooth structure as well.  The precise definition is given in 
the next section.  

\noindent The main theorem is:

\begin{thm} Let  $\Gamma \subset SL(2, {\bf Z})$ be a 
subgroup containing two hyperbolic elements.  We assume
that the (four) eigenvectors of these two
elements are pairwise linearly independent. Let $0 < \alpha < 1$.
Then the $C^{1 + \alpha}$ \T \ space ${\rm Teich}^{1 + \alpha}
(F: \Gamma \times {\bf T}^{n} 
\rightarrow {\bf T}^{n})$ is trivial, where $F$ is the standard
action. 
\end{thm}

\noindent {\bf Remark.} Any non-cyclic subgroup satisfies the hypotheses 
of the Theorem.

\bigskip

\noindent An immediate corollary of Theorem 1 is:

\begin{thm} Let $0 < \alpha < 1$. The $C^{1 + \alpha}$ \T \ space 
${\rm Teich}^{1 + \alpha}
(F: SL(2,{\bf Z}) \times {\bf T}^{2} \rightarrow  
{\bf T}^2)$ is trivial, where $F$ is the standard action.
\end{thm}

\noindent {\bf Acknowledgements.}  I thank Dennis Sullivan and Michael 
Lyubich for helpful conversations, and Geoffrey Mess and Steve Hurder
for making me aware of this problem.

\bigskip \bigskip

\noindent {\bf 2} \ \ \ {\bf The \T \ space of an Anosov action}

\bigskip

\noindent We consider $\Gamma$, a locally compact, 
second countable group, and a
$C^{r}$ left action 
\[
F: \Gamma \times M  \rightarrow M
\]
of $\Gamma$ on a $C^{r}$ manifold $M$, $1 \leq r \leq \omega$.
The action is called Anosov if at least one element $\gamma \in
\Gamma$ acts by an Anosov diffeomorphism.  That is, there is a continous
splitting of the tangent bundle $TM = E^{s} \bigoplus E^{u}$, invariant
under $D\gamma$, and such that vectors in $E^{s}$ are
exponentially contracted, and vectors in 
$E^{u}$ are exponentially expanded, by iterates of $D\gamma$.

The \T \ space of $F$ 
is defined as follows.  Consider triples $(h, N, G)$ where $N$ is a 
$C^{r}$ manifold, $h: M \rightarrow N$ is a homeomorphism, and 
\[
G:  \Gamma  \times N \rightarrow N
\]
is a $C^{r}$ action conjugate to $F$ by $h$.  That is, 
$G(x,\gamma) = h(F( h^{-1}(x), \gamma))$, for every $\gamma \in 
\Gamma$.  We assume in addition that if $F(  \gamma, \cdot)$ is 
Anosov, then $G( \gamma, \cdot)$ is also Anosov.  Such a triple
is called a {\em marked Anosov action modeled on $F$}.  Two triples
$(h_{1}, N_{1}, G_{1})$ and
$(h_{2}, N_{2}, G_{2})$ are {\em equivalent} if the homeomorphism
$s : N_{1} \rightarrow N_{2}$ defined by $s \circ h_{1} = s_{2}$ is
a $C^{1}$ diffeomorphism.

\bigskip \bigskip

\noindent {\bf 3} \ \ \ {\bf The main argument}

\bigskip

\noindent The subbundles $E^{s}$ and $E^{u}$ in the definition of an Anosov 
map are integrable, and the corresponding foliations are called 
the stable, respectively unstable, foliations of the map.
The diffeomorphism is called codimension-one if either the stable
or the unstable foliation is codimension-one. If the 
diffeomorphism is $C^{1 + \alpha}$, i.e. the derivative is \Ho \
continuous with exponent $0 < \alpha < 1$,  then the codimension-one foliation 
has $C^{1 + \beta}$ transverse regularity for some $0 < \beta < 1$.\cite{Ho}, \cite{Mn}
(A codimension-$k$ Anosov foliation
is transversely absolutely continuous, and the holonomy maps have 
\Ho \ Jacobian. \cite{An})
An Anosov diffeomorphism of ${\bf T}^{2}$ has simultaneous foliation 
charts
\[
\phi: D^1 \times D^1 \rightarrow U \subset {\bf T}^{2}
\]
where $D^1$ is the one dimensional disk. 
The intersection of a leaf of the unstable foliation ${\cal W}^{u}$ with
the neighborhood $U$ is a union of horizontals $\phi(D^{1} \times {y})$,
and that of a leaf of the stable foliation ${\cal W}^{s}$ is a union of
verticals $\phi({x} \times D^{1})$.  Moreover, we can choose $\phi$ to 
be smooth along, say $x_{0} \times D^1$ and $D^{1} \times y_{0}$, for 
some $(x_{0}, y_{0}) \in D^{1} \times D^{1}$.
A simple but important 
observation is the following:  since both foliations have 
$C^{1 + \beta}$
transverse regularity, these charts belong to the $C^{1 + \beta}$ smooth
structure on ${\bf T}^{2}$ preserved by the diffeomorphism.
Therefore the smooth structure is determined (up to 
$C^{1 + \beta}$ equivalence) by the pair of transverse
smooth structures.

We give the main argument of the proof of Theorem 1.  The lemmas used 
here are proved in the next section.  Let $\gamma_{1}$, $\gamma_2 \in
SL(2, {\bf Z})$  be hyperbolic elements such that 
the eigenvectors $\lbrace v_{1}^{s}, v_{1}^{u}, v_{2}^{s},
v_{2}^{u} \rbrace$ are pairwise linearly independent.  Here $v_{i}^s$,
respectively $v_{i}^{u}$, is the contracting (or stable), respectively
expanding (or unstable), eigenvector of $\gamma_i$, for $i = 1,2$.  Let 
\[
F: \Gamma \times {\bf T}^{2}  \rightarrow {\bf T}^{2}
\]
be the standard linear action.  Define
\[
f_{i}( \cdot ) = F( \gamma_{i}, \cdot)
\]
for $i = 1,2$.  
The lines in ${\bf R}^2$ parallel to $v_{i}^{s}$ 
project to ${\bf T}^{2} = {\bf R}^{2} / {\bf Z}^{2}$ to give the
stable foliation ${\cal W}^{s,F}_{i}$, for $i = 1,2$.   
The  unstable foliation ${\cal W}^{u,F}_{i}$ is obtained analogously.

Consider a segment $\tau$ contained in a leaf $W \in {\cal W}_{1}^{u,F}$.
Then $\tau$ is a transversal to both ${\cal W}_{1}^{s,F}$ and ${\cal W}_
{2}^{s,F}$.  There is a {\em locally defined} holonomy from $\tau$ to itself,
``generated'' by the pair of foliations as follows.  Slide $\tau$ a 
small distance along the leaves of ${\cal W}_{1}^{s,F}$, remaining 
transverse to both stable foliations.  Then slide $\tau$ up the leaves
of ${\cal W}_{2}^{s,F}$, returning to the original leaf $W$ of 
${\cal W}_{1}^{u,F}$. This is all done inside a single foliation chart for 
${\cal W}_{1}^{u,F}$.  The resulting motion is simply rigid translation in
the leaf $W$.  In other words, the translation group of $W$, restricted
to small translations defined on the segment $\tau$, can be canonically
factored into a composition of holonomy along ${\cal W}_{1}^{s,F}$,
followed by holonomy along ${\cal W}_{2}^{s,F}$.

\begin{figure}
\vspace{3.5in}
\caption{Translations along $W$ are compositions of holonomy.}
\end{figure}

Let 
\[
G: \Gamma \times {\bf T}^{2} \rightarrow {\bf T}^{2}
\]
define a point in the  $C^{1 + \alpha}$ \T \ space of the 
standard action of $\Gamma$ on
${\bf T}^{2}$.  In other words, $G$ is a $C^{1 + \alpha}$ action, and there
is a homeomorphism $h: {\bf T}^{2} \rightarrow
{\bf T}^{2}$ such that 
\[
G(\gamma,x) = h ( F( \gamma,h^{-1}(x)))
\]
for all $\gamma \in \Gamma$.  Let 
\[
g_{i}( \cdot) = G( \gamma_{i}, \cdot)
\]
for $i = 1,2$. Then $g_{1}$ and $g_{2}$ are Anosov.
Let ${\cal W}_{i}^{s,G}$ and ${\cal W}_{i}^{u,G}$ be the stable
and unstable foliations respectively of $g_{i}$.  We have
$h({\cal W}_{i}^{s,F}) = {\cal W}_{i}^{s,G}$ and $h({\cal W}_{i}
^{u,F}) = {\cal W}_{i}^{u,G}$.  We consider $W^{\prime} = h(W)$,
and $\tau^{\prime} = h(\tau)$.  The homeomorphism $h$ conjugates
the small translations on $\tau \subset W$ to an action on 
$\tau^{\prime} \subset W^{\prime}$.  
The main claim is that this action, denoted 
\[
S: [0,\epsilon] \times \tau^{\prime} \rightarrow \tau^{\prime},
\]
is smooth. More precisely,  
\begin{lem}
$\frac{\partial S}{\partial y} (t,y)$ is continuous.  
\end{lem}
The following
rigidity property of the translation pseudogroup  of the line
implies that $h$ must be $C^1$ along $\tau$.
\begin{prop}
Let $T: {\bf R} \times {\bf R} \rightarrow {\bf R}$ be the rigid
translations, namely $T(t,x) = x + t$.  Let $S: {\bf R} \times {\bf R}
\rightarrow {\bf R}$ be an action conjugate to $T$ by a homeomorphism
$h: {\bf R} \rightarrow {\bf R}$.  So $S(t,y) = h ( T (t, h^{-1}(y)))$.
If $S$ is smooth in a weak sense, namely if 
$\frac{\partial S}{\partial y} (t,y)$ is continuous, then in fact
$h$ is a $C^1$ diffeomorphism.
\end{prop}
To prove Lemma 1, we need the following important fact.

\begin{lem}
The foliations \WonesG, \WoneuG, \WtwosG, \WtwouG are pairwise transverse.
\end{lem}
Of course we know \WonesG and \WoneuG are transverse, and similarly
\WtwosG and \WtwouG are transverse.  
But one doesn't know a priori that the homeomorphism $h$ did not
introduce tangencies between, say, \WoneuG \ and \WtwouG.  

\bigskip

\noindent {\bf Remark.}  Lemma 2 is automatically satisfied for actions
{\em sufficiently near} to the standard action, since the transversality
property is stable.  The fact that this can be {\em proved} without
any nearness assumption is the ingredient that globalizes the rigidity
result.

\bigskip

\noindent {\bf Proof of Lemma 1, assuming Lemma 2.}
The conjugacy $h$ preserves the foliations.  Thus 
we note that the action $S$ on $\tau^{\prime}$ can be factored
into holonomy along \WonesG \ followed by holonomy along \WtwosG, simply
by applying $h$ to the corresponding decomposition for the rigid
translations on $\tau$.  The foliations \WonesG \ and  \WtwosG \ are
transversely $C^{1 + \beta}$.  The leaves of \WoneuG \ form a 
$C^{1 + \beta}$ family of transversals  to \WonesG \ and  \WtwosG.
The claimed regularity of $S$ follows.

\begin{flushright}
$\Box$
\end{flushright}

The argument up to now shows that the transverse smooth structure to the 
unstable foliation \WoneuG \ is $C^1$ equivalent by the conjugacy
$h$ to the tranverse smooth structure to ${\cal W}^{u,F}_{1}$, the
corresponding foliation for the standard action.
Similarly we see that the transverse smooth structure to the stable
foliation \WonesG \ is $C^1$ equivalent by $h$ to the linear one. 
But the simultaneous foliation charts are charts in the $C^{1 + \beta}$
smooth structure, hence 
we conclude that $h$ is $C^1$.

\noindent {\bf Remark.}  It can be shown that 
$h$ is in fact as smooth as the mapping $G$.\cite{L}

\bigskip \bigskip

\noindent {\bf 4} \ \ \ {\bf The translation pseudogroup of the 
real line is rigid}

\bigskip

\noindent We prove Proposition 1.
Let $T: {\bf R} \times {\bf R} \rightarrow {\bf R}$  be the rigid
translations, $T(t,x) = x + t$.  Let $S: {\bf R} \times {\bf R} 
\rightarrow {\bf R}$ be conjugate to $T$ by the homeomorphism $h$.
We assume that $\frac{\partial S}{\partial y}(t,y)$ is continuous.
We will show that $S$ is $C^1$ conjugate to a rigid
translation $L: {\bf R} \times {\bf R} \rightarrow {\bf R}$,
$L(t,z) = z + \alpha t$, where $\alpha$ is a fixed real number.  Since
$L$ is affinely conjugate to $T$, we will conclude that $h$ is $C^1$.

\bigskip

Pick $y_{0} \in {\bf R}$. Let $t(y)$ be defined by $S(t(y),y) = y_{0}$.
Define $g: {\bf R} \rightarrow {\bf R}$ 
\[
g(y) = \int_{y_0}^{y} \frac{\partial S}{\partial y^{\prime}}
(t(y^{\prime}), y^{\prime})dy^{\prime}.
\]
Conjugate the action of $S$ by the homeomorphism $g$ to obtain an
action $L$, 
\[
L(t,z) = g(S(t,g^{-1}(z))).
\]
One sees that $\frac{\partial L}{\partial z}(t,z) = 1$.  Hence 
$L(t,y) = y + k(t)$, for some function $k(t)$.  But
$L( t + s,y) = L(s,L(y,t)) = L(s,y + k(t)) = y + k(t) + k(s)$, 
so $k(t) = \alpha t$ for some $\alpha$.

\bigskip \bigskip

\noindent {\bf 5} \ \ \ {\bf The foliations remain transverse}

\bigskip

\noindent In this section we prove Lemma 2.
There is an equivalence relation naturally associated to an 
Anosov diffeomorphism, defined by the pair of foliations
${\cal W}^s$ and ${\cal W}^u$.  Namely, $x$ is equivalent to $y$ if $x \in 
W^{s}(y) \cap W^{u}(y)$.  This equivalence relation is generated
by a pseudogroup of local homeomorphisms, 
defined in simultaneous foliation charts as the
product of holonomy along ${\cal W}^{u}$ by holonomy along ${\cal W}^{s}$.

\begin{lem}
Suppose \WoneuG \ is tangent to \WtwosG \ at $z$.  The \WoneuG
\ is tangent to \WtwosG \ at every $z^{\prime} \in 
W^{u,G}_{1}(z) \cap W^{s,G}_{1}(z)$.
\end{lem}
Assuming this, we can prove Lemma 2.
The set of $z^{\prime} \in W^{u,G}_{1}(z) \cap W^{s,G}_{1}(z)$
is dense in ${\bf T}^{2}$, and the distributions $E_{1}^{u,G}$
and $E_{2}^{s,G}$ are continuous.  But the foliations 
\WoneuG \ and \WtwosG \ do not coincide, since they are 
homeomorphic images of the transverse foliations ${\cal W}^{u,F}_{1}$  
and ${\cal W}^{s,F}_{2}$.  Therefore there can be no tangencies, and 
Lemma 2 is proved.

\bigskip

\noindent {\bf Proof of Lemma 3.}
Let $z^{\prime} \in W^{u,G}_{1}(z) \cap W^{s,G}_{1}(z)$.
We denote by $W^{u,G}_{1}(z, \epsilon)$ a small neighborhood 
of $z$ in the leaf of \WoneuG \ containing it.
Let 
\[
hol_{s}: W^{u,G}_{1}(z, \epsilon) \rightarrow 
W^{u,G}_{1}(z^{\prime}, \epsilon)
\]
be the map defined by the holonomy of $W^{s,G}_{1}$.  Similarly,
let 
\[
hol_{u}: W^{s,G}_{1}(z, \epsilon) \rightarrow
W^{s,G}_{1}(z^{\prime}, \epsilon)
\]
be the holonomy of \WoneuG.  We can represent $W^{s,G}_{2}(z, \epsilon)$
as the graph of a map
\[
\theta_{z}: W^{u,G}_{1}(z, \epsilon) \rightarrow W^{s,G}_{1}(z,\epsilon).
\]
Similarly we can represent $W^{s,G}_{2}(z^{\prime},\epsilon)$ as the 
graph of a map 
\[
\theta_{z^{\prime}}: W^{u,G}_{1}(z^{\prime},\epsilon) \rightarrow
W^{s,G}_{1}(z^{\prime}, \epsilon).
\]
We observe that
\[
\theta_{z^{\prime}} = hol_{s} \circ \theta_{z} \circ hol_{u}^{-1}.
\]
To see this, consider the corresponding graphs and holonomy for the 
linear map.  The corresponding equation holds there.  This is a 
consequence of a non-obvious feature of the foliations in the linear case:
the local pseudogroup of homeomorphisms that generates the Anosov 
equivalence relation described above {\em preserves the 
contracting (or expanding) foliation of any other hyperbolic element}.
Since the equation is defined by conjugacy invariant objects,
it must hold for the non-standard action as well.  Since
$hol_{s}$ and $hol_{u}$ are $C^{1 + \beta}$ diffeomorphisms, we
see that a tangency at $z$ forces a tangency at $z^{\prime}$.

\begin{flushright}
$\Box$
\end{flushright}

\bibliographystyle{alpha}

\bibliography{ref}



\begin{thebibliography}{MM87b}

\bibitem[An]{An}
D.V. Anosov.
\newblock Geodesic flows on closed {Riemannian} manifolds with negative
  curvature.
\newblock {\em Proc. Inst. Stek.}, 90:1 -- 235, 1967.

\bibitem[C]{C}
E.~Cawley.
\newblock The {Teichm\"{u}ller} space of an {Anosov} diffeomorphism of ${\bf
  T}^2$.
\newblock Inst. for Math. Sci. at Stony Brook, preprint 1991/9.

\bibitem[L]{L}
R.~de~la Llave.
\newblock Invariants for smooth conjugacy of hyperbolic dynamical systems {II}.
\newblock {\em Communications in Mathematical Physics}, 109:369 -- 378, 1987.

\bibitem[HKLZ]{HKLZ}
S.~Hurder, A.~Katok, J.~Lewis, and R.~Zimmer.
\newblock Rigidity for {Cartan} actions of higher rank lattices.
\newblock preprint.

\bibitem[Ho]{Ho}
E.~Hopf.
\newblock Statistik der geod\"{a}tischen {Linien} in {Mannigfaltigkeiten} negativer
  {Kr\"{u}mmung}.
\newblock {\em Ber. Verh. Sochs. Akad. Wiss. Leipzig}, 91:261 -- 304, 1939.

\bibitem[Hu]{Hu}
S.~Hurder.
\newblock Rigidity for {Anosov} actions of higher rank lattices.
\newblock preprint.

\bibitem[KL]{KL}
A.~Katok and J.~Lewis.
\newblock Local rigidity of certain groups of toral automorphisms.
\newblock preprint.

\bibitem[MM1]{MM1}
J.M. Marco and R.~Moriyon.
\newblock Invariants for smooth conjugacy of hyperbolic dynamical systems {I}.
\newblock {\em Communications in Mathematical Physics}, 109:681 -- 689, 1987.

\bibitem[MM2]{MM2}
J.M. Marco and R.~Moriyon.
\newblock Invariants for smooth conjugacy of hyperbolic dynamical systems
  {III}.
\newblock {\em Communications in Mathematical Physics}, 112:317 -- 333, 1987.

\bibitem[Mn]{Mn}
R. Ma\~{n}\'e.
\newblock {\em Ergodic Theory and Differentiable Dynamics}.
\newblock Springer-Verlag, New York, 1983.

\bibitem[Ta]{Ta}
F.~Tangerman.
\newblock short note.

\end{thebibliography}

\end{document}